\documentclass[12pt,reqno]{amsart}
\usepackage{amssymb}
\usepackage{}
\usepackage{amsfonts}
\usepackage{mathrsfs}
\usepackage{amsmath,amsthm,amssymb,amsfonts,amscd}
\usepackage{mathrsfs}
\usepackage{bbding}
\usepackage{graphicx,latexsym}

\theoremstyle{plain}
\newtheorem{theorem}{Theorem}[section]

\newtheorem{lemma}{Lemma}[section]

\newtheorem{conjecture}{Conjecture}[section]

\theoremstyle{definition}

\def \N {{\mathbb N}}
\def \Z {{\mathbb Z}}

\def\bg{\bigg}
\def\({\bg(}
\def\){\bg)}

\def\le {\leqslant}
\def\ge {\geqslant}

\begin{document}
%\hbox{Preprint}
\medskip

\title{On a conjecture of Sun involving powers of three}

\author{Quan-Hui Yang}
\address{School of Mathematics and Statistics \\ Nanjing University of Information, Science and Technology \\ Nanjing 210044 \\Reople's Republic of China}
\email{yangquanhui01@163.com}
\author{Lilu Zhao}
\address{School of Mathematics \\ Shandong University \\ Jinan  250100 \\Reople's Republic of China}
\email{zhaolilu@sdu.edu.cn}

\begin{abstract} Given a positive integer $n\ge 2$, let $D(n)$ denote the smallest positive integer $m$ such that $a^3+a(1\le a\le n)$ are pairwise distinct modulo $m^2$. A conjecture of Z.-W. Sun states that $D(n)=3^k$, where $3^k$ is the least power of $3$ no less than $\sqrt{n}$. The purpose of this paper is to confirm this conjecture.
\end{abstract}

\maketitle

\section{Introduction}
\setcounter{lemma}{0}
\setcounter{theorem}{0}
\setcounter{corollary}{0}
\setcounter{remark}{0}
\setcounter{equation}{0}

\medskip

Let $f(x)\in \Z[x]$ be a polynomial with all $f(a)(a\in \Z^{+})$ pairwise distinct. Define the determinant $\Delta_f(n)$ to be the smallest positive integer $m$ such that $f(a)(1\le a\le n)$ are pairwise distinct modulo $m$.

In 1985 Arnold, Benkoski and McCabe \cite{ABM} studied $\Delta_f(n)$ with $f(x)=x^2$, and it was shown in \cite{ABM} that for $n>4$, $\Delta_f(n)$ is the smallest positive integer $m\ge 2n$ such that $m$ is $p$ or $2p$ with $p$ an odd prime. Sun \cite{S13} determined $\Delta_f(n)$ for $f(x)=2x(x-1)$ and $f(x)=x(x-1)$ respectively. For example, Sun \cite{S13} proved that if $f(x)=2x(x-1)$, then $\Delta_{f}(n)$ is the least prime number greater than $2n-2$. For the study of $\Delta_f(n)$ with higher degree polynomials $f$, one may refer to \cite{BSW,Moree,Zieve}. Sun proposed many interesting conjectures related to determinants, and one may refer to \cite{S13} (see also Section 4 of Chapter 6 in \cite{S2021}) for details.

In this paper, we consider $D(n)$, which denotes the smallest positive integer $m$ such that $a^3+a(1\le a\le n)$ are pairwise distinct modulo $m^2$. In 2013 Z.-W. Sun made the following conjecture (see Conjecture 6.76 in \cite{S2021}).
\begin{conjecture} \label{conjecture1}For all $n\ge 2$, one has
$D(n)=3^k$, where
\begin{align}\label{definek}k=\min\{j\in \Z^{+}:\ 3^j\ge \sqrt{n}\}.\end{align}
\end{conjecture}
The main result in this paper is to confirm the above conjecture.
\begin{theorem}\label{theorem1}Conjecture \ref{conjecture1} is true.\end{theorem}
A closely related problem is to deal with $\Delta_f(n)$ with $f(x)=x^3+x$, and the method in this paper can be also used to study it. We shall consider it elsewhere.

\section{Preparations}

\setcounter{lemma}{0}
\setcounter{theorem}{0}
\setcounter{corollary}{0}
\setcounter{remark}{0}
\setcounter{equation}{0}

We first point out that $\sqrt{n}\le D(n)\le 3^k$, where $k$ is given in \eqref{definek} throughout this paper. In fact, we have the following.
\begin{lemma}\label{lemma1}Let $n\ge 2$. Then $a^3+a(1\le a\le n)$ are pairwise distinct modulo $3^{2k}$.\end{lemma}
\begin{proof} Suppose that $1\le a<b\le n$. Note that $b^3+b-a^3-a=(b-a)(a^2+ab+b^2+1)$. Since $3\nmid (a^2+ab+b^2+1)$ and $1\le b-a<n\le 3^{2k}$, one has
$b^3+b\not\equiv a^3+a\pmod{3^{2k}}$. This completes the proof.
\end{proof}

By Lemma \ref{lemma1}, in order to establish Theorem \ref{theorem1}, it suffices to prove the following lemma.
\begin{lemma}\label{lemma2}Let $n\ge 2$. Suppose that
\begin{align}\label{inequality} \sqrt{n}<m<3^k<3\sqrt{n}.\end{align}
Then there exist  $1\le a<b\le n$ such that $b^3+b\equiv a^3+a\pmod{m^2}$.\end{lemma}
In order to prove Lemma \ref{lemma2}, we shall consider the following six cases.

(i) $m=\delta p^r$, where $1\le \delta \le 3$, $p\ge 5$ is a prime and $r\in \Z^{+}$.

(ii) $m=2^r$, where $r\in \Z^{+}$.

(iii) $m=2^r3^{s}$, where $r,s\in \Z^{+}$.

(iv) $m=2^r3^{s}5$, where $r,s\in \N=\Z^{+}\cup\{0\}$.

(v) $m=2^r3^{s}7$ or $m=2^r\cdot 3^{s}\cdot 5\cdot 7$, where $r,s\in \N$.

(vi) $m=\delta p^{r}$, where $\delta \ge 4$, $p\ge 5$ is a prime, $r\in \Z^{+}$, $p\nmid \delta$ and $p^r\ge 11$.

\noindent It is clear that $m$ in \eqref{inequality} must satisfy (at least) one of the above six cases.

\bigskip

\section{The Cases (ii)-(vi)}

\setcounter{lemma}{0}
\setcounter{theorem}{0}
\setcounter{corollary}{0}
\setcounter{remark}{0}
\setcounter{equation}{0}

\medskip

The purpose of this section is to deal with cases (ii)-(vi). Throughout this paper, we assume that \eqref{inequality} holds.
\begin{lemma}\label{lemma31}Suppose that $m=\delta p^{r}$, where $\delta \ge 4$, $p\ge 5$ is a prime, $r\in \Z^{+}$ and $p\nmid \delta$. If $p^{2r}+\delta^2p\le n$, then there exist  $1\le a<b\le n$ such that $b^3+b\equiv a^3+a\pmod{m^2}$.\end{lemma}
\begin{proof}We consider $1\le a\le p^{2r}$ and $b=a+\delta^2c$ with $1\le c\le p$. Note that $b\le p^{2r}+\delta^2p\le n$. It suffices to find $a,c$ such that
$$a^2+a(a+\delta^2c)+(a+\delta^2c)^2+1\equiv 0\pmod{p^{2r}}.$$
This is equivalent to
$$(6a+3\delta^2c)^2\equiv -3\delta^4c^2-12\pmod{p^{2r}}.$$
It is well known that
$$\Big|\sum_{1\le j\le p}\Big(\frac{-3\delta^4j^2-12}{p}\Big)\Big|=1,$$
where $\big(\frac{\cdot}{p}\big)$ denotes the Legendre symbol. We conclude that there exists $1\le c\le p$ such that $-3\delta^4c^2-12$ is a quadratic residue modulo $p$. Now it is easy to deduce that there exists $1\le a\le p^{2r}$ such that $(6a+3\delta^2c)^2\equiv -3\delta^4c^2-12\pmod{p^{2r}}$. This completes the proof.
\end{proof}

\begin{lemma}[Case (vi)]\label{lemma32}Let $n\ge 2000$. Suppose that $m=\delta p^{r}$, where $\delta \ge 4$, $p\ge 5$ is a prime, $r\in \Z^{+}$, $p\nmid \delta$ and $p^r\ge 11$. Then there exist  $1\le a<b\le n$ such that $b^3+b\equiv a^3+a\pmod{m^2}$.\end{lemma}
\begin{proof}In view of Lemma \ref{lemma31}, we need to verify $p^{2r}+\delta^2p\le n$. By \eqref{inequality}, $n> \frac{\delta^2p^{2r}}{9}$. It suffices to prove  $p^{2r}+\delta^2p\le \frac{\delta^2p^{2r}}{9}$, which is equivalent to $(\delta^2-9)(p^{2r-1}-9)\ge 81$. Since $\delta^2-9\ge 7$ and $p^{2r-1}\ge p^{r}\ge 11$, we have
$$(\delta^2-9)(p^{2r-1}-9)\ge \max\{2(\delta^2-9),\ 7(p^{2r-1}-9)\}.$$
Note that $\delta^2 (p^{2r-1})^2\ge \delta^2p^{2r}=m^2>n\ge 2000$. Thus, we have either $\delta^2\ge 50$ or $p^{2r-1}\ge 21$, and it follows that $(\delta^2-9)(p^{2r-1}-9)\ge 81$. This completes the proof.
\end{proof}

\begin{lemma}[Case (v)]\label{lemma33}Let $n\ge 14$. Suppose that $m=2^r3^{s}7$ or $m=2^r\cdot3^{s}\cdot 5\cdot 7$, where $r,s\in \N$. Then there exist  $1\le a<b\le n$ such that $b^3+b\equiv a^3+a\pmod{m^2}$.\end{lemma}
\begin{proof}We first consider the case $r\ge 1$. We write $m^2=14t$. Note that $14|t$. We choose $a=3$, $b=3+t$ and it is easy to verify $b^3+b\equiv a^3+a\pmod{m^2}$. Since $m^2=14t<9n$, we have
$3+t<3+\frac{9}{14}n\le n$.

Now we assume that $r=0$. Thus $m=3^s7$ or $m=3^s\cdot 5\cdot 7$. We write $m^2=7t$ and choose $a=3$, $b=3+t$. Then again it is easy to verify $b^3+b\equiv a^3+a\pmod{m^2}$. We need to confirm $3+t\le n$. If $m=3^s7$, then by \eqref{inequality} we have $s=k-2$ and $3+t=3+3^{2s}7=3+3^{2k-4}7<3+\frac{7}{9}n<n$.  If $m=3^s\cdot 5\cdot 7$, then by \eqref{inequality} we have $s=k-4$ and $3+t=3+3^{2s}5^27=3+3^{2k-8}5^27<3+\frac{175}{3^6}n<n$.
This completes the proof.
\end{proof}

\begin{lemma}\label{lemma34}Suppose that $m=2^{r}t$, where $r\in \Z^{+}$, $t$ is an odd number. If $2^{2r}+t^2\le n$, then there exist  $1\le a<b\le n$ such that $b^3+b\equiv a^3+a\pmod{m^2}$.\end{lemma}
\begin{proof}We consider $1\le a\le 2^{2r}$ and $b=a+t^2$. Note that $b\le 2^{2r}+t^2\le n$. It suffices to find $a$ such that
$$a^2+a(a+t^2)+(a+t^2)^2+1\equiv 0\pmod{2^{2r}}.$$
This is equivalent to
$$3a^2+3at^2+t^4+1\equiv 0\pmod{2^{2r}}.$$
In fact, we can prove by induction that for any $j\in \Z^{+}$, there exists $1\le a\le 2^{j}$ such that
\begin{align}\label{induction}3a^2+3at^2+t^4+1\equiv 0\pmod{2^{j}}.\end{align}
When $j=1$, we choose $a=1$. Suppose that \eqref{induction} holds. We consider $a'=a+2^jc$. Then
\begin{align*}3a'^2+3a't^2+t^4+1\equiv 3a^2+3at^2+t^4+1+2^jc\pmod{2^{j+1}}.\end{align*}
We can find $c\in \{0,1\}$ such that $3a'^2+3a't^2+t^4+1\equiv 0\pmod{2^{j+1}}$.
This completes the proof.
\end{proof}

\begin{lemma}[Case (iv)]\label{lemma35}Let $n\ge 57$. Suppose that $m=2^r3^{s}5$, where $r,s\in \N$. Then there exist  $1\le a<b\le n$ such that $b^3+b\equiv a^3+a\pmod{m^2}$.\end{lemma}
\begin{proof}We first consider the case $r\ge 2$. In view of Lemma \ref{lemma34}, we only need to verify
$2^{2r}+t^2\le n$, where $t=3^s5$. This follows from
$2^{2r}+t^2\le \frac{2^{2r}t^2}{9}$, because $n>\frac{m^2}{9}=\frac{2^{2r}t^2}{9}$. The inequality $2^{2r}+t^2\le \frac{2^{2r}t^2}{9}$ holds due to $2^{2r}\ge 16$ and $t^2\ge 25$.

Now we assume that $r=0$ or $r=1$. Then $m=3^s5$ or $m=2\cdot 3^s \cdot 5$. We write $m=\delta p$ with $p=5$. In view of Lemma \ref{lemma31}, we only need to verify $25+5\delta^2\le n$. If $m=3^s5$, then by \eqref{inequality}, we have $s=k-2$ and $25+5t=25+5\cdot 3^{2s}=25+5\cdot 3^{2k-4}<25+\frac{5}{9}n\le n$. If $m=2\cdot 3^s\cdot 5$, then by \eqref{inequality}, we have $s=k-3$ and $25+5t=25+20\cdot 3^{2s}=25+20\cdot 3^{2k-6}<25+\frac{20}{81}n\le n$.
This completes the proof.
\end{proof}

\begin{lemma}[Case (iii)]\label{lemma36}Let $n\ge 144$. Suppose that $m=2^r3^{s}$, where $r,s\in \Z^{+}$. Then there exist  $1\le a<b\le n$ such that $b^3+b\equiv a^3+a\pmod{m^2}$.\end{lemma}
\begin{proof}We first consider the case $r\ge 2$ and $s\ge 2$. It is easy to verify
$2^{2r}+3^{2s}\le n$, and the desired conclusion follows from Lemma \ref{lemma34}.

Next we consider the case $r=1$ and $s\ge 2$.  By \eqref{inequality}, we have $s=k-1$ and $1+3^{2s}=1+ 3^{2k-2}<1+\frac{n}{4}<n$. The desired conclusion follows by choosing $a=1$ and $b=1+3^{2s}$, since $a^2+ab+b^2+1\equiv 0\pmod{4}$.

Now we assume that $s=1$. Then $m=2^r3$ with $r\ge 2$. Note that $a^2+a(a+9)+(a+9)^2+1$ is equal to $112=2^4\cdot 7$ if $a=1$. Similarly to \eqref{induction}, we can prove that for any $j\ge 4$, there exists $1\le a\le 2^{j}-15$ such that
\begin{align*}a^2+a(a+9)+(a+9)^2+1\equiv 0\pmod{2^{j}}.\end{align*}
In particular, there exists $1\le a\le 2^{2r}-15$ such that
$a^2+a(a+9)+(a+9)^2+1\equiv 0\pmod{2^{2r}}$. The desired conclusion follows by choosing $b=a+9$.
This completes the proof.
\end{proof}

\begin{lemma}[Case (ii)]\label{lemma37}Let $n\ge 64$. Suppose that $m=2^r$, where $r\in \Z^{+}$. Then there exist $1\le a<b\le n$ such that $b^3+b\equiv a^3+a\pmod{m^2}$.\end{lemma}
\begin{proof}Note that $(a+4)^3+(a+4)-a^3-a=4(3(a+2)^2+5)$. It suffices to find $1\le a\le 2^{2r-4}-3$ such that $3(a+2)^2+5\equiv 0\pmod{2^{2r-2}}$.

Note that $3\cdot 1^2+5\equiv 0\pmod{8}$. We can prove by induction that for any $j\ge 3$, there exists $x$ such that \begin{align}\label{induction2}3x^2+5\equiv 0\pmod{2^j}.\end{align}
From \eqref{induction2}, for $j\ge 3$, we can deduce that $3(x+2^{j-1}c)^2+5\equiv 3x^2+5+2^jc \pmod{2^{j+1}}$. Thus, there exists $x'$ such that $3x'^2+5\equiv 0\pmod{2^{j+1}}$. In particular, there exists $x$ such that $3x^2+5\equiv 0\pmod{2^{2r-2}}$. Since $x$ is odd, we can assume that $1\le x\le 2^{2r-3}-1$.

Let $y=2^{2r-3}-x$. Then we deduce that $3y^2+5=3(2^{2r-3}-x)^2+5\equiv 3x^2+5\equiv 0\pmod{2^{2r-2}}$. Therefore, we can further assume that $1\le x\le 2^{2r-4}-1$. Since $n\ge 64$, we have $r\ge 4$. Then we obtain  $3\le x\le 2^{2r-4}-1$.  On choosing $a=x-2$, we obtain $3(a+2)^2+5\equiv 0\pmod{2^{2r-2}}$. Note that $a+4\le 2^{2r-4}+1<\frac{9}{16}n+1\le n$.
This completes the proof.
\end{proof}

\bigskip

\section{The Cases (i)}

\setcounter{lemma}{0}
\setcounter{theorem}{0}
\setcounter{corollary}{0}
\setcounter{remark}{0}
\setcounter{equation}{0}

\medskip

In this section, we deal with the case (i). The first result is the following.
\begin{lemma}\label{lemma41}Let $n\ge 100$. Suppose that $m=\delta p^{r}$, where $1\le \delta \le 3$, $p\equiv 1\pmod{3}$ is a prime, $r\in \Z^{+}$. Then there exist  $1\le a<b\le n$ such that $b^3+b\equiv a^3+a\pmod{m^2}$.\end{lemma}
\begin{proof}We consider $1\le a\le p^{r}$ and $b=a+\delta^2p^r$. Note that
$$b\le p^{r}+\delta^2p^r=(1+\delta^2)p^r\le(1+\delta^2)\frac{3\sqrt{n}}{\delta}\le 10\sqrt{n}\le n.$$
 It suffices to prove that there exists $1\le a\le p^r$ such that
$$a^2+a(a+\delta^2p^r)+(a+\delta^2p^r)^2+1\equiv 0\pmod{p^{r}}.$$
This is equivalent to
$$3a^2+1\equiv 0\pmod{p^{r}}.$$
The desired conclusion now follows due to the condition $p\equiv 1\pmod{3}$. This completes the proof.
\end{proof}

From now on, we assume throughout this section that $m=\delta p^{r}$, where $1\le \delta \le 3$, $p\equiv 2\pmod{3}$ is an odd prime, $r\in \Z^{+}$.
Note that $p\equiv 2\pmod{3}$ is equivalent to $\big(\frac{-3}{p}\big)=-1$.

We shall deal with the remaining case by using analytic number theory method. Throughout this section, we use the notation
$$e(\alpha)=e^{2\pi i\alpha}.$$

We introduce
\begin{align}X:=X_p=\begin{cases}[\frac{p^2}{9}]p^{2r-2}, \ & \textrm{ if } p=5,
\\ [\frac{p}{9}]p^{2r-1}, \ & \textrm{ if } p\ge 11.\end{cases}\end{align}
We write
\begin{align}\label{definerho}\rho:=\rho_p=\begin{cases}2r-2, \ & \textrm{ if } p=5,
\\ 2r-1, \ & \textrm{ if } p\ge 11.\end{cases}\end{align}

We aim to find $1\le a\not=b\le \frac{n}{\delta^2}$ such that $\delta^4(a^2+ab+b^2)+1\equiv 0\pmod{p^{2r}}$. Then on choosing $a'=\delta^2 a$, $b'=\delta^2 b$, we obtain $a'^3+a'\equiv b'^3+b'\pmod{m^2}$. By \eqref{inequality}, we have $X<\frac{n}{\delta^2}$.

Let
\begin{align}\label{definef}f(a,b)=\delta^4(a^2+ab+b^2).\end{align}Now we introduce
\begin{align}\label{defineN}\mathcal{N}=\sum_{\substack{1\le a,b\le X \\ f(a,b)+1\equiv 0\pmod{p^{2r}}}}1.\end{align}
Since $p\equiv 2\pmod{3}$ is an odd prime, one has $f(a,a)+1\not\equiv 0\pmod{p^{2r}}$. In particular,
\begin{align}\label{Nwithnotequal}\mathcal{N}=\sum_{\substack{1\le a\not=b\le X \\ f(a,b)+1\equiv 0\pmod{p^{2r}}}}1.\end{align}
Therefore, the main objective is to prove $\mathcal{N}>0$.

For $j\ge 1$, we define
\begin{align}\label{defineTj}T_j=\sum_{\substack{1\le c\le p^j \\ (c,p)=1}}\sum_{\substack{1\le a,b\le X }}e\big(\frac{cf(a,b)+c}{p^j}\big).\end{align}

\begin{lemma}\label{lemma42}Let $\mathcal{N}$ and $T_j$ be given in \eqref{defineN} and \eqref{defineTj} respectively. We have
\begin{align}\mathcal{N}=\frac{X^2}{p^{2r}}+\frac{1}{p^{2r}}\sum_{j=1}^{2r}T_j.\end{align}
\end{lemma}
\begin{proof}We make use of the identity
\begin{align}\label{identity}\frac{1}{q}\sum_{c=1}^{q}e\big(\frac{cg}{q}\big)=\begin{cases}
1, \ & \textrm{ if } q|g
\\ 0, \ & \textrm{ if } q\nmid g\end{cases}\end{align}
to obtain
\begin{align*}\mathcal{N}=\frac{1}{p^{2r}}\sum_{\substack{1\le c\le p^{2r} }}\sum_{\substack{1\le a,b\le X }}e\big(\frac{cf(a,b)+c}{p^{2r}}\big).\end{align*}
We change variable by taking $c=p^{2r-j}c'$ with $0\le j\le 2r$, $1\le c'\le p^{j}$ and $(c',p)=1$ to deduce that
\begin{align*}\mathcal{N}=\frac{1}{p^{2r}}\sum_{j=0}^{2r}\sum_{\substack{1\le c'\le p^{j} \\ (c',p)=1 }}\sum_{\substack{1\le a,b\le X }}e\big(\frac{c'f(a,b)+c'}{p^{j}}\big)=\frac{X^2}{p^{2r}}+\frac{X^2}{p^{2r}}\sum_{j=1}^{2r}T_j.\end{align*}
We are done.
\end{proof}

We define
\begin{align}\label{defineSj}S_j:=S_j(x,y)=\sum_{\substack{1\le c\le p^j \\ (c,p)=1}}\sum_{\substack{1\le a,b\le p^j }}e\big(\frac{cf(a,b)+ax+by+c}{p^j}\big).\end{align}

\begin{lemma}\label{lemma43}Let $S_j$ be given in \eqref{defineSj}. We have
\begin{align}\label{Sj}S_j=p^j\Big(\frac{-3}{p^j}\Big)\sum_{\substack{1\le c\le p^j \\ (c,p)=1}}e\big(\frac{\overline{3c\delta^4}(-x^2+xy-y^2)+c}{p^j}\big),\end{align}
where the notation $\overline{d}$ in \eqref{Sj} means $d\cdot \overline{d}\equiv 1\pmod{p^j}$.
\end{lemma}
\begin{proof}Write \begin{align}\label{defineRj}R_j=\sum_{\substack{1\le a,b\le p^j }}e\big(\frac{cf(a,b)+ax+by}{p^j}\big).\end{align}
The desired conclusion follows from
\begin{align}\label{Rj}R_j=p^j\big(\frac{-3}{p^j}\big)e\big(\frac{\overline{3c\delta^4}(-x^2+xy-y^2)}{p^j}\big),\end{align}
It is not hard to verify that
\begin{align*}&c\delta^4(a^2+ab+b^2)+ax+by
\\ \equiv &  c\delta^4(a+\overline{2}b+\overline{2c\delta^4}x)^2+3c\delta^4(\overline{2}b+\overline{3c\delta^4}(y-\overline{2}x))^2
+\overline{3c\delta^4}(-x^2+xy-y^2)\pmod{p^j}.
\end{align*}
We conclude from above that
\begin{align*}R_j=\Big(\sum_{\substack{1\le a\le p^j }}e\big(\frac{ca^2}{p^j}\big)\Big)\Big(\sum_{\substack{1\le b\le p^j }}e\big(\frac{3cb^2}{p^j}\big)\Big)e\big(\frac{\overline{3c\delta^4}(-x^2+xy-y^2)}{p^j}\big).\end{align*}
For the Gauss sum, it is well known that (see Chapter 7 of \cite{Hua})
$$\sum_{\substack{1\le a\le p^j }}e\big(\frac{ca^2}{p^j}\big)=p^{\frac{j}{2}}\Big(\frac{c}{p^j}\Big)\epsilon_{p^j},$$
where $\epsilon_{p^j}$ satisfies $\epsilon_{p^j}^2=\big(\frac{-1}{p^j}\big)$. This proves \eqref{Rj} and we are done.
\end{proof}

\begin{lemma}\label{lemma44}Let $\rho$ be given in \eqref{definerho}. For $1\le j\le \rho$,  we have
\begin{align}\label{Tj-1}T_j=X^2p^{-j}\big(\frac{-3}{p^j}\big)\mu(p^j),\end{align}
where $\mu(\cdot)$ is the M\"obius function.
\end{lemma}
\begin{proof}If $1\le j\le \rho$, then $p^j|X$ and we have
\begin{align}T_j=X^2p^{-2j}\sum_{\substack{1\le c\le p^j \\ (c,p)=1}}\sum_{\substack{1\le a,b\le p^j }}e\big(\frac{cf(a,b)+c}{p^j}\big)=X^2p^{-2j}S_j(0,0).\end{align}
By Lemma \ref{lemma43}, we have $S_j(0,0)=p^{j}\big(\frac{-3}{p^j}\big)\mu(p^j)$ and this yields \eqref{Tj-1}.
We are done.\end{proof}

\begin{lemma}\label{lemma45}Let $\rho$ be given in \eqref{definerho}. We have
\begin{align}\mathcal{N}=\frac{X^2}{p^{2r}}(1+\frac{1}{p})+\frac{1}{p^{2r}}\sum_{j=\rho+1}^{2r}T_j.\end{align}
\end{lemma}
\begin{proof}The desired conclusion follows from Lemma \ref{lemma42} and Lemma \ref{lemma44}.\end{proof}

Note that for $\rho+1\le j\le 2r$, the summations over $a$ and $b$ in \eqref{defineTj} are incomplete summations (modulo $p^j$).

\begin{lemma}\label{lemma46}For $\rho+1\le j\le 2r$, we have
\begin{align}\label{Tj-2}T_j=\frac{1}{p^{2j}}\sum_{\substack{1\le x,y\le p^j}}S_j(x,y)\sum_{\substack{1\le a',b'\le X}}e\big(-\frac{a'x+b'y}{p^j}\big).\end{align}
\end{lemma}
\begin{proof}This follows from \eqref{identity} by considering the summations over $x$ and $y$ in \eqref{Tj-2}.\end{proof}

\begin{lemma}\label{lemma47}Let $u\in \Z$. We have
\begin{align}\label{boundKl}\Big|\sum_{\substack{1\le c\le p^j \\ (c,p)=1}}e\big(\frac{\overline{c}u+c}{p^j}\big)\Big|\le 2p^{\frac{j}{2}}.\end{align}
In particular, we have
\begin{align}\label{Sjinequality}|S_j|\le 2p^{\frac{3j}{2}}.\end{align}
\end{lemma}
\begin{proof}Write
\begin{align}K(p^j;u)=\sum_{\substack{1\le c\le p^j \\ (c,p)=1}}e\big(\frac{\overline{c}u+c}{p^j}\big).\end{align}
If $p^j|u$, then $K(p^j;u)$ is a Ramanujan sum and $K(p^j;u)=\mu(p^j)$. If $p\nmid u$, then $K(p^j;u)$ is a Kloosterman sum and by Corollary 4.4 in \cite{Iwaniec} we have $|K(p^j;u)|\le 2p^{\frac{j}{2}}$.

Now we assume that $p^t\| u$ with $1\le t\le j-1$. By changing variables $c=xp^{j-1}+y$, we obtain
\begin{align}K(p^j;u)=\sum_{\substack{1\le y\le p^j \\ (y,p)=1}}\sum_{1\le x\le p}e\big(\frac{\overline{y}u+xp^{j-1}+y}{p^j}\big)=0.\end{align}
The proof of \eqref{boundKl} is complete. Note that \eqref{Sjinequality} follows from Lemma \ref{lemma43} and \eqref{boundKl} immediately. We are done.
 \end{proof}

\begin{lemma}\label{lemma48}For $\rho+1\le j\le 2r$, we have
\begin{align}\label{bound}\sum_{\substack{1\le x\le p^j}}\Big|\sum_{\substack{1\le a\le X}}e\big(\frac{ax}{p^j}\big)\Big|\le p^{j}(2+\ln p^j).\end{align}
\end{lemma}
\begin{proof} On writing $Y=\frac{p^j-1}{2}$, we have
\begin{align*}\sum_{\substack{1\le x\le p^j}}\Big|\sum_{\substack{1\le a\le X}}e\big(\frac{ax}{p^j}\big)\Big|=&X+2\sum_{\substack{1\le x\le Y}}\Big|\sum_{\substack{1\le a\le X}}e\big(\frac{ax}{p^j}\big)\Big|
\\ = &X+2\sum_{\substack{1\le x\le Y}}\Big|\frac{1-e\big(\frac{xX}{p^j}\big)}{1-e\big(\frac{x}{p^j}\big)}\Big|
\\ \le & X+4\sum_{\substack{1\le x\le Y}}\Big|\frac{1}{1-e\big(\frac{x}{p^j}\big)}\Big|.\end{align*}
For $0<\theta<1/2$, we have $|1-e(\theta)|=2\sin (\pi \theta)>4\theta$. Therefore,
\begin{align*}\sum_{\substack{1\le x\le p^j}}\Big|\sum_{\substack{1\le a\le X}}e\big(\frac{ax}{p^j}\big)\Big|
\le & X+\sum_{\substack{1\le x\le Y}}\frac{p^j}{x}\le X+p^j(1+\ln Y)\le p^{j}(2+\ln p^j).\end{align*}We are done.
\end{proof}

\begin{lemma}\label{lemma49}For $\rho+1\le j\le 2r$, we have
\begin{align}\label{boundTj}|T_j|\le 2p^{\frac{3j}{2}}(2+\ln p^j)^2.\end{align}

\end{lemma}
\begin{proof}By Lemmas \ref{lemma46}-\ref{lemma48}, we deduce that
\begin{align*}|T_j|\le &\frac{1}{p^{2j}}\sum_{\substack{1\le x,y\le p^j}}|S_j(x,y)|\cdot \Big|\sum_{\substack{1\le a',b'\le X}}e\big(-\frac{a'x+b'y}{p^j}\big)\Big|
\\ \le & 2p^{\frac{3j}{2}} \cdot \frac{1}{p^{2j}}\sum_{\substack{1\le x,y\le p^j}} \Big|\sum_{\substack{1\le a',b'\le X}}e\big(-\frac{a'x+b'y}{p^j}\big)\Big|
\\= & 2p^{\frac{3j}{2}} \cdot \frac{1}{p^{2j}}\Big(\sum_{\substack{1\le x\le p^j}}\Big|\sum_{\substack{1\le a\le X}}e\big(\frac{ax}{p^j}\big)\Big|\Big)^{2}
\\ \le & 2p^{\frac{3j}{2}}(2+\ln p^j)^2.\end{align*}
We are done.
\end{proof}

\begin{lemma}\label{lemma410}We have
\begin{align}\label{finalineq}\Big|\mathcal{N}-\frac{X^2}{p^{2r}}\big(1+\frac{1}{p}\big)\Big|\le
\begin{cases}2p^{r}(2+\ln p^{2r})^2(1+\frac{1}{p\sqrt{p}}),\ & \textrm{ if } p=5,
\\ 2p^{r}(2+\ln p^{2r})^2, \ & \textrm{ if } p\ge 11. \end{cases}\end{align}
\end{lemma}
\begin{proof}The desired conclusion follows from Lemma \ref{lemma45} and Lemma \ref{lemma49}.\end{proof}

\begin{lemma}\label{lemma411}Let $n\ge 2\cdot 10^{12}$. Then we have
\begin{align}\mathcal{N}>0.\end{align}
\end{lemma}
\begin{proof}For $p=5$, by \eqref{finalineq}, we need to prove
$$\frac{[\frac{p^2}{9}]^2p^{4r-4}(1+\frac{1}{p})}{p^{2r}}>8p^r(1+\ln p^r)^2(1+\frac{1}{p\sqrt{p}}),$$
and this follows from
\begin{align}\label{check1}p^{r}>2\cdot 5^4(1+\ln p^r)^2.\end{align}

For $p\ge 23$, we $[\frac{p}{9}]>\frac{p}{9}-1\ge \frac{1}{17}p$ and thus
\begin{align}\label{check2}[\frac{p}{9}]^2p^{-2}\ge 17^{-2}.\end{align}
One can check that \eqref{check2} holds for $11\le p\le 19$ as well.

For $p\ge 11$, by \eqref{finalineq}, we need to prove
$$\frac{[\frac{p}{9}]^2p^{4r-2}}{p^{2r}}>8p^r(1+\ln p^r)^2,$$
and by \eqref{check2}, this follows from
\begin{align}\label{check3}p^{r}>8\cdot 17^2(1+\ln p^r)^2.\end{align}
On writing $q=\sqrt{p^r}$, our task is to prove $q>34\sqrt{2}(1+2\ln q)$.

Let $f(x)=x-34\sqrt{2}(1+2\ln x)$. Then
$f'(x)=1-68\sqrt{2}\cdot\frac{1}{x}$ and $f$ is increasing when $x> 68\sqrt{2}$. Note that $f(680)>0$.
By \eqref{inequality}, $q=\sqrt{p^r}>\sqrt{\frac{\sqrt{n}}{\delta}}\ge \sqrt{\frac{\sqrt{n}}{3}}>680$. This establishes \eqref{check1} and \eqref{check3}.
The proof is complete.
\end{proof}

\begin{lemma}[Case (i)]\label{lemma412}Let $n\ge 2\cdot 10^{12}$. Suppose that $m=\delta p^{r}$, where $1\le \delta \le 3$, $p\ge 5$ is a prime, $r\in \Z^{+}$. Then there exist $1\le a<b\le n$ such that $b^3+b\equiv a^3+a\pmod{m^2}$.\end{lemma}
\begin{proof}If $p\equiv 1\pmod{3}$, then the desired conclusion follows from Lemmas \ref{lemma41}. If $p\equiv 2\pmod{3}$, then the desired conclusion follows by combining \eqref{Nwithnotequal} and Lemma \ref{lemma411}. \end{proof}
\noindent {\it Proof of Lemma \ref{lemma2}.} We assume that $n\ge 2\cdot 10^{12}$ (and the conclusion for $n<2\cdot 10^{12}$ can be checked directly). Lemma \ref{lemma2} follows from Lemma \ref{lemma32}, Lemma \ref{lemma33}, Lemma \ref{lemma35}, Lemma \ref{lemma36}, Lemma \ref{lemma37} and Lemma \ref{lemma412}. According to the remark between Lemma \ref{lemma1} and Lemma \ref{lemma2}, we also complete the proof of Theorem \ref{theorem1}.

\end{document}